\documentclass[11pt]{article}

\usepackage{lscape,tikz}
\usepackage{epsfig,cancel,ulem, amsthm,amssymb}
\definecolor{darkblue}{rgb}{0.2,0.2,0.71}
\usepackage{todonotes}

\usepackage{color,soul}
\usetikzlibrary{decorations.markings,arrows}

\usepackage{wrapfig}
\usepackage{color}
\usepackage{graphicx,framed}
\usepackage[colorlinks=true, pdfstartview=FitV, linkcolor=darkblue, citecolor=darkblue, urlcolor=darkblue]{hyperref}
\definecolor{shadecolor}{rgb}{0.95, 0.95, 0.86}
\definecolor{darkgreen}{rgb}{0.2, 0.5,  0}

\def\&{\vspace{-5pt}&}

\newtheorem{theorem}{Theorem}[section]
\newtheorem{example}[theorem]{Example}
\newtheorem{exercise}[theorem]{Exercise}

\newtheorem{lemma}[theorem]{Lemma}
\newtheorem{remark}[theorem]{Remark}

\newtheorem{proposition}[theorem]{Proposition}
\newtheorem{corollary}[theorem]{Corollary}
\newtheorem{definition}[theorem]{Definition}

\def\bt{\begin{theorem}}
\def\et{\end{theorem}}
\def\bc{\begin{corollary}}
\def\ec{\end{corollary}}
\def\bx{\begin{example}}
\def\ex{\end{example}}
\def\bxr{\begin{exercise}\small}
\def\exr{\end{exercise}}
\def\bl{\begin{lemma}}
\def\el{\end{lemma}}
\def\bd{\begin{definition}}
\def\ed{\end{definition}}
\def\bp{\begin{proposition}}
\def\ep{\end{proposition}}

\def\br{\begin{remark}}
\def\er{\end{remark}}

\def\be{\begin{equation}}
\def\ee{\end{equation}}

\def\&{\hspace{-15pt}&}
\def\bea{\begin{eqnarray}}
\def\eea{\end{eqnarray}}

 \newcommand{\C}{\mathbb{C}}
 \newcommand{\W}{\mathrm{W}}

\makeatletter
\@addtoreset{equation}{section}
\makeatother

\begin{document}
\title{Dicyclic groups and Frobenius manifolds}
\author{Yassir Dinar \and Zainab Al-Maamari }
\date{}

\maketitle

\begin{abstract}
The orbits space of an irreducible representation of a finite group is a variety whose coordinate ring is finitely generated by homogeneous invariant polynomials. Boris Dubrovin showed that  the orbits spaces of the reflection groups acquire the structure of polynomial Frobenius manifolds. Dubrovin's method to construct examples of Frobenius manifolds on orbits spaces  was carried for other linear representations of discrete groups which have in common that the coordinate rings of the the orbits spaces are  polynomial rings. In this article, we show that the orbits space of an irreducible representation of a Dicyclic group acquire two structures of Frobenius manifolds. The coordinate ring of this orbits space is not a polynomial ring.
\end{abstract}

\maketitle
\tableofcontents

\section{Introduction}

The notion of Frobenius manifold was introduced by Boris Dubrovin as a geometric realization of a potential $\mathbb F$ satisfying a  system of partial differential equations known in topological fields theory as WDVV equations \cite{DuRev}.  Beside topological fields theory,  Frobenius manifolds appear in many fields like invariant theory, integrable systems, quantum cohomology and singularity theory. This article contributes to the  relation between Frobenius manifolds and  invariant theory.

Let $\W$ be a finite group of linear transformations acting  on a complex vector space $V$ of dimention $r$.  Then the orbits space $M=V/\W$ of this group is a variety whose coordinate ring is the ring of invariant polynomials $\C [V]^\W$. The ring $\C [V]^\W$ is finitely generated by homogeneous polynomials.  If $f_1,f_2,\ldots,f_m$ is a set of such  generators then $m\geq r$ and the relation between them is called syzygies. The set of generators are not unique, nor are their degrees \cite{Neusel},\cite{derksen}.

An element $w\in \W$ is called a reflection if it fixes a subspace of $V$ of codimention one pointwise. The group $\W$ is called a complex reflection group  if it is generated by reflections.   Then  Shephard-Todd-Chevalley theorem states that $\W$ is a reflection group if and only if the invariant ring $\C [V]^\W$ is a polynomial ring \cite{Neusel}, i.e. it is generated by $r$ algebraically independent homogeneous polynomials (so there are no syzygies).  Furthermore, when $\W$ is a reflection group, the degrees of such a set generators of  $\C [V]^\W$ are uniquely specified by the group and we   refer to them as   the degrees of $\W$.

Let us assume $\W$ is a Shephard group, i.e. a symmetry group of a regular complex polytope. Then $\W$ is a reflection group. Let $f_1,f_2,\ldots,f_r$ be a set of algebraically independent homogeneous generators of  $\C [V]^\W$.  We assume that    degree $f_i$ is less than or equal degree $f_j$  when $i<j$. Then the inverse of the Hessian of $f_1$ defines a flat metric $(\cdot,\cdot)_2$ on $T^*M$ \cite{orlik}.  There is another  flat metric $(\cdot,\cdot)_1$ on $T^*M$, which was  studied initially by K. Saito (\cite{saito1}, \cite{saito2}),  defined  as the Lie derivative of $(\cdot,\cdot)_2$ along the vector field $e=\partial_{f_r}$. The two metrics form  what is called a flat pencil of metrics (more details are given below). Dubrovin used the properties of this flat pencil of metrics to construct polynomial Frobenius manifolds \cite{DubAlmost} (see \cite{DCG} and \cite{zuo1} for the case of Coxeter groups). This article is about applying Dubrovin's method for other finite linear groups than Shephard groups.

Dubrovin's method to construct Frobenius manifolds, through finding flat pencils of metrics on orbits spaces,  was carried out for infinite linear groups like extended affine Weyl groups \cite{dubext}, \cite{dubext1}, Jacobi groups \cite{bertola} and recently a new extension of affine Weyl groups \cite{zuo2}. They all have in common that the  invariant rings  are  polynomial rings. Moreover, even when  considering a generalization of Frobenius manifold structure  on orbits spaces, many results was  obtained under the assumption that the invariant ring is a polynomial ring \cite{arsie}. Then it is natural question to ask about applying Dubrovin's method on orbits spaces of  finite non-reflection groups.

In this article we apply Dubrovin's method and construct Frobenius manifolds on orbits spaces of Dicyclic groups. The resulting  Frobenius manifolds can be obtained by using an ad-hoc procedure, but it is fascinating to find them on orbits spaces of some group. Precisely,
we will show that the orbits space of the Dicyclic group of order $4n$ is endowed with two structure of  Frobenius manifolds which up to scaling has the following potential
\be
\mathbb F(z_1,z_2)=z_1^{k}+\frac{1}{2} z_2^2 z_1, ~k={3-d\over 1-d}
\ee
where $d$ is $\frac{2 + \sqrt{3}n}{\sqrt{3}n}$ or $\frac{2 - \sqrt{3}n}{\sqrt{3}n}$.

To make the article as self-contained as possible, we review in next  section the definition of Frobenius manifold and its relation with the theory of flat pencils of metrics. In the last section we obtain the promised Frobenius manifolds by direct calculations.

\section{Preliminaries}
\subsection{\textbf{Frobenius manifolds}}
A Frobenius algebra is a commutative associative algebra with unity $e$ and an invariant nondegenerate bilinear form $<\cdot,\cdot>$.
A \textbf{Frobenius manifold} is a manifold $M$ with
a smooth structure of a Frobenius algebra on the tangent space $T_tM$ at
any point $t \in M $ with certain compatibility conditions \cite{DuRev}. Globally, we  require the metric $<\cdot,\cdot>$ to be flat and the unity vector field $e$ to be covariantly constant with respect to it. In the flat  coordinates $(t^1,...,t^r)$ where $e={\partial\over \partial t^{r}}$ the compatibility conditions imply  that there exists a function $\mathbb{F}(t^1,...,t^r)$ such that
\[ \eta_{ij}=<\partial_{t^i},\partial_{t^j}>=  \partial_{t^{r}}
\partial_{t^i}
\partial_{t^j} \mathbb{F}(t)\]
and the structure constants of the Frobenius algebra are given by
\[ C_{ij}^k=\sum_p\Omega_1^{kp}  \partial_{t^p}\partial_{t^i}\partial_{t^j} \mathbb{F}(t)\]
where $\Omega_1^{ij}$ denote the inverse of the matrix $\eta_{ij}$.
In this work, we consider Frobenius manifolds where the quasihomogeneity condition takes the form
\begin{equation}
\sum_{i=1}^r d_i t^i \partial_{t^i} \mathbb{F}(t) = \left(3-d \right) \mathbb{F}(t);~~~~d_{r}=1.
\end{equation}
 This condition defines {\bf the degrees} $d_i$ and {\bf the charge} $d$ of  the Frobenius structure.
The  associativity of the Frobenius
algebra implies that the potential $\mathbb{F}(t)$ satisfies a system of partial differential equations which appears in topological field theory
and is called Witten-Dijkgraaf-Verlinde-Verlinde (WDVV)   equations:
\begin{equation} \label{frob}
\sum_{k,p} \partial_{t^i}
\partial_{t^j}
\partial_{t^k} \mathbb{F}(t)~ \Omega_1^{kp} ~\partial_{t^p}
\partial_{t^q}
\partial_{t^n} \mathbb{F}(t) =\sum_{k,p} \partial_{t^n}
\partial_{t^j}
\partial_{t^k} \mathbb{F}(t) ~\Omega_1^{kp}~\partial_{t^p}
\partial_{t^q}
\partial_{t^i} \mathbb{F}(t),~~ \forall i,j,q,n.
  \end{equation}
 Detailed information about Frobenius manifolds and related topics can be found in $\cite{DuRev}$.

\subsection{Flat pencil of metrics and Frobenius manifolds}
In this section we review the relation between the geometry  of flat pencil of metrics and Frobenius manifolds. See  \cite{DFP} for details.

Let $M$ be a smooth manifold of dimension $r$. A symmetric bilinear form $(\cdot,\cdot)$ on $T^*M$ is called a \textbf{contravariant
metric} if it is invertible on an open dense subset $M_0 \subset M$. In local coordinates $(u^1, . . . , u^r)$, if we set
\be
\Omega^{ij}(u)=(du^i, du^j);~ i, j = 1, . . . , r.
\ee
Then  the inverse matrix $\Omega_{ij}(u)$ of $ \Omega^{ij}(u)$ determines a metric $<\cdot,\cdot >$ on $TM_0$. We define the \textbf{contravariant Christoffel symbols} $\Gamma^{ij}_k$ of $(\cdot,\cdot )$ by $
\Gamma^{ij}_k:=-\sum_s{\Omega}^{is} \Gamma_{sk}^j$
where $\Gamma_{sk}^j$ are the  Christoffel symbols of $<\cdot,\cdot>$. We say the metric $(\cdot,\cdot)$ is flat if  $<\cdot,\cdot>$ is flat.

Let  $(\cdot,\cdot)_1$ and $(\cdot,\cdot)_2$ be two contravariant flat metrics on $M$ and denote their  Christoffel symbols by $\Gamma_{1;k}^{ij}(u)$ and $\Gamma_{2;k}^{ij}(u)$ respectively.   We say $(\cdot,\cdot)_1$ and $(\cdot,\cdot)_2$ form a  \textbf{flat pencil of metrics} if  $(\cdot,\cdot)_\lambda:=(\cdot,\cdot)_1+\lambda (\cdot,\cdot)_2$ defines a flat metric on $T^*M$ for a generic $\lambda$ and
 its   Christoffel symbols are given by $\Gamma_{\lambda ;k}^{ij}(u)=\Gamma_{2;k}^{ij}(u)+\lambda \Gamma_{1;k}^{ij}(u)$.

Let $(\cdot,\cdot)_1$ and $(\cdot,\cdot)_2$  be two contravariant metrics on $M$ and denote their matrices by  ${\Omega}_1^{ij}(u)$ and ${\Omega}_2^{ij}(u)$, respectively, in some coordinates $(u^1,\ldots,u^r)$. Suppose that they form  a flat pencil of metrics. This flat pencil of metrics is called  \textbf{ quasihomogeneous of  degree} $d$ if there exists a
function $\tau$ on $M$ such that the vector fields
\begin{eqnarray}\label{fpm}  E&:=& \nabla_2 \tau, ~~E^i
=\sum_s{\Omega}_2^{is}\partial_s\tau
\\\nonumber  e&:=&\nabla_1 \tau, ~~e^i
= \sum_s {\Omega}_1^{is}\partial_s\tau  \end{eqnarray} satisfy the following
relations
\[ [e,E]=e,~~ \textrm{Lie}_E (~,~)_2 =(d-1)(~,~)_2,~~ \textrm{Lie}_e (~,~)_2 =
(~,~)_1,~~ \textrm{Lie}_e(~,~)_1
=0.\]
Here $\textrm{Lie}_X$ denote the Lie derivative  along a given  vector field $X$.  In addition, the  quasihomogeneous flat pencil of metrics is called \textbf{regular} if  the
(1,1)-tensor $
  R_i^j = {d-1\over 2}\delta_i^j + {\nabla_1}_i
E^j$
is  nondegenerate on $M$.

The following theorem due to Dubrovin  gives a connection between the geometry of Frobenius manifolds and flat pencils of metrics.
\bt \cite{DFP}\label{dub flat pencil}
A  quasihomogeneous regular  flat pencil of metrics of degree $d$ on a manifold $M$ defines a Frobenius structure on $M$ of  charge $d$.
 \et

 Let us assume the flat pencil of metrics on $M$ is  regular quasihomogeneous   of degree $d$. Let $(t^1,\ldots,t^r)$ be flat coordinates of $(\cdot,\cdot)_1$ where  $\tau=t^1$, $e=\partial_{t^r}$ and $E=\sum_i d_i t^i {\partial_{t^i}}$. Let $\eta_{ij}$ denote the inverse of $\Omega_1^{ij}(t)$. Then it turns out that the potential $\mathbb F(t^1,\ldots,t^r)$ is obtained from the equations
 \begin{eqnarray}
\frac{\partial^2 \mathbb F}{\partial t^i \partial t ^j}&=&\sum_{k,l}\frac{1}{d-1+d_k+d_l}\eta_{ik}\eta_{jl}\Omega_2^{kl}(t)\\\nonumber
\frac{\partial \mathbb F}{\partial t^i}&=&\frac{1}{(3-d-d_i)}\sum d_k t^k \frac{\partial^2 \mathbb F}{\partial t^k \partial t^i}\\\nonumber
\mathbb F(t)&=&\frac{1}{3-d} \sum d_k t^k \frac{\partial \mathbb F}{\partial t^k}\nonumber
\end{eqnarray}

It is well known that from a Frobenius manifold we always have a flat pencil of metrics but it does not necessarily satisfy the regularity condition \cite{DFP}.

\section{Dicyclic groups}

Let $n$ be a natural number greater that 1 and  $\W$ be the matrix  group generated by
\be
\sigma:=\left(\begin{array}{cc}
     \xi& 0  \\
0     & \xi^{-1}
\end{array}\right),~\alpha :=\left(\begin{array}{cc}
    0 & 1 \\
    -1 & 0
\end{array}\right)
\ee
where $\xi$ is a primitive $2n$-th root of unity. Then $\sigma$ and $\alpha$ satisfy the relations
\be
\sigma^{2n}=1, \alpha^2=\sigma^n, \alpha^{-1} \sigma \alpha=\sigma^{-1}.
 \ee
Thus $\W$ is isomorphic to the dicyclic group  of order $4n$. The invariant ring of $\W$ is   generated by the following homogeneous polynomials \cite{Neusel}
\be u_1=x_1^2x_2^2,~ u_2= x_1^{2n}+x_2^{2n},~u_3=x_1x_2(x_1^{2n}-x_2^{2n})
\ee
subject to  the relation
\begin{equation} \label{defeq}
u_3^2-u_1u_2^2+4 u_1^{n+1}=0.
\end{equation}

The orbits space $M$ of $\W$ is a variety  isomorphic to  the  hypersurface  $T$ defined as the  zero set of equation (\ref{defeq}) in $\mathbb{C}^3$. Consider  equation (\ref{defeq}) as a quadratic equation in $u_3$. Then any point $p$ out of the discriminant locus has small neighbourhood  $U_p$ where   $u_1$ and $u_2$ act as  coordinates. In what follows we assume that we fix such  open set $U\subset V$  with    coordinates $(u_1,u_2)$ and we omit the subscript $p$.

Let $h$ be the Hessain matrix of $u_1$, i.e. $h_{ij}={\partial^2 u_1\over \partial x_i\partial x_j}$ and let $h^{-1}$ denotes its inverse. Then, by direct calculations,  $h^{-1}$  defines a flat contravariant metric $(\cdot,\cdot)_2$ on $U$.  This metric, in the coordinates $u_1$ and $u_2$, is given as follows
\be
(\cdot,\cdot)_2=
\left(
\begin{array}{cc}
 \frac{4 }{3}u_1 & \frac{2n}{3}  u_2 \\
 \frac{2n}{3}  u_2 & -\frac{2 n^2}{3u_1} (  u_2^2-6 u_1^n)
\end{array}
\right)
;~~  (dt_i,dt_j)_2=\sum_{k,l=1}^2{\partial u_i\over \partial x_k}{\partial u_j\over \partial x_l}{h_{kl}^{-1}}.
\ee

 Let $e$ be a vector field of the form $f(u_1)\partial_{u_2}$, where $f(u_1)$ is any  smooth function. Then, by direct calculations,  the  Lie derivative $(\cdot,\cdot)_1$  of $(\cdot,\cdot)_2$ along  $e$ forms  with $(\cdot,\cdot)_2$ a flat pencil of metrics.    This metric  takes the value  \be
(\cdot,\cdot)_1=\left(
\begin{array}{cc}
 0 & \frac{2}{3} \left(n f-2 u_1 f'\right) \\
 \frac{2}{3} \left(n f-2 u_1 f'\right) & -\frac{4 }{3 u_1}\left(n^2u_2f  +n u_1 u_2 f'\right)
\end{array}
\right).
\ee
The guess for the vector field to take this from was inspired by \cite{almaamari}. In order to get a quasihomogeneous flat pencil of metrics, we need the Lie derivative of $(\cdot,\cdot)_1$ with respect to  $e$ to equal zero. This condition leads to the following differential  equation for $f(u)$
\be
2 n u_1 f f'-2 u_1^2 (f')^2+n^2 f^2=0
\ee
which has two independent solutions
\be
f_+=u_1^{\frac{n}{2} \left(1+\sqrt{3}\right)} \textrm{ and } f_-=u_1^{\frac{n}{2} \left(1-\sqrt{3}\right)}.
\ee
Let us assume  $e=f_+\partial_{u_2}=u_1^{\frac{n}{2} \left(1+\sqrt{3}\right)}\partial_{u_2}$. Then
\be
(dt_i,dt_j)_1=\left(
\begin{array}{cc}
 0 & -\frac{2 n}{\sqrt{3}} u_1^{\frac{1}{2} \left(1+\sqrt{3}\right) n} \\
 -\frac{2 n}{\sqrt{3}} u_1^{\frac{1}{2} \left(1+\sqrt{3}\right) n} & -\frac{2}{3} \left(3+\sqrt{3}\right) n^2 u_1^{\frac{1}{2} \left(1+\sqrt{3}\right) n-1} u_2
\end{array}
\right).
\ee
It turns out that  the two metrics $(\cdot,\cdot)_2$ and $(\cdot,\cdot)_1$  form a quasihomogeneous flat pencil of metrics with degree
\be
d=\frac{2 + \sqrt{3}n}{\sqrt{3}n}.
\ee
In the notations of equations (\ref{fpm}), we have  $\tau=-{\sqrt{3}\over 2n}u_1$ and
\be
E=-\frac{2 }{\sqrt{3} n}u_1\partial_{u_1}-\frac{1}{\sqrt{3}}u_2\partial_{u_2}.
\ee
This flat pencil of metrics is also regular since the $(1,1)$-tensor $R_i^j$ equals  the nondegenerate  matrix
\be
\left(\begin{array}{cc}
     -\frac{1}{\sqrt{3}n}& 0  \\
     0 & \frac{1-n}{\sqrt{3}n}
\end{array}\right).
\ee
Flat coordinates for $(\cdot,\cdot)_1$ are obtained by setting
\be t_1=-\frac{\sqrt{3}}{2 n} u_1,~~t_2=u_2 u_1^{-\frac{1}{2} \left(1+\sqrt{3}\right)
   n}.
\ee
In these coordinates we get
\be
(\cdot,\cdot)_2=\left(\begin{array}{cc}
 -\frac{2}{\sqrt{3} n} t_1 & t_2 \\
 t_2 & 2^{1-\sqrt{3} n} 3^{\frac{1}{2} \left(\sqrt{3} n+1\right)} n^2 \left(-n
   t_1\right){}^{-\sqrt{3} n-1}
\end{array}
\right),~~(\cdot,\cdot)_1=\left(
\begin{array}{cc}
 0 & 1 \\
 1 & 0
\end{array}
\right)
\ee
The  potential $\mathbb F_+$ of the corresponding Frobenius manifold is
\be\label{posPotential}
{\mathbb F_+}=\frac{2^{-\sqrt{3} n}
   3^{\frac{1}{2}
   \left(\sqrt{3} n+1\right)}
   \left(-n
   t_1\right){}^{1-\sqrt{3}
   n}}{3 n^2-1}+\frac{1}{2}
   t_1 t_2^2
\ee

Let us  take $e=f_-\partial_{u^2}=u_1^{\frac{n}{2} \left(1-\sqrt{3}\right)}\partial_{u_2}$. Then  similar to above,  we get a regular quasihomogenous flat pencil of metrics of degree
\be
d=\frac{2 - \sqrt{3}n}{\sqrt{3}n}
\ee
with $\tau={\sqrt{3}\over 2n}u_1$. The resulting potential will be
\be \label{negPotential}
\mathbb F_-=
\frac{2^{\sqrt{3} n}
   3^{\frac{1}{2}-\frac{\sqrt{3} n}{2}}
   \left(n t_1\right){}^{\sqrt{3} n+1}}{3
   n^2-1}+\frac{1}{2} t_1 t_2^2.
   \ee

We  repeat the calculation by taking $(u_1,u_3)$ as coordinates instead of $(u_1,u_2)$.  It turns out that even though the middle steps may differ in values, the resulting Frobenius manifolds are exactly the same as those given by the potentials (\ref{posPotential}),(\ref{negPotential}).

We observe that Dubrovin computed by ad hoc procedure  all possible potentials of 2-dimensional Frobenius manifolds \cite{DuRev}. The potentials found in this article, after scaling,   are listed by Dubrovin in the form
\be
\mathbb F(z_1,z_2)=z_1^{k}+\frac{1}{2} z_2^2 z_1, ~k={3-d\over 1-d}
\ee
where $d$  is $\frac{2 + \sqrt{3}n}{\sqrt{3}n}$ or $\frac{2 - \sqrt{3}n}{\sqrt{3}n}$. However, finding it by using the method of a flat pencil of metrics on an orbits space of a finite group that is not a reflection group is a surprising result.

The result reported in this article is a part of work in progress to apply Dubrovin's method on orbits spaces of finite groups to find new interesting examples of Frobenius manifolds. In future publications, we will consider irreducible representations of Coxeter groups which are not  reflection representations  \cite{almaamari}.

\noindent{\bf Acknowledgments.}

 The authors thank
Hans-Christian Herbig  for stimulating discussions.    This work is funded by the internal grant of Sultan Qaboos University (IG/SCI/DOMS/19/08). The
authors like to thank anonymous reviewers  for their comments and suggestions.

\noindent Yassir Dinar \\
\noindent dinar@squ.edu.om \\

\noindent Zainab Al-Maamari\\
\noindent z.mamari2@gmail.com\\

\noindent Depatment of Mathematics\\
\noindent College of Science\\
\noindent Sultan Qaboos University\\
\noindent Muscat, Oman

\end{document}